\documentclass[11pt,reqno]{amsart}

\numberwithin{equation}{section}
\usepackage{amssymb,rotating,amsmath,amsfonts,amsthm,color,graphics,epsfig,bbm,bm}

\newcommand{\calR}{\mathcal{R}}

\newcommand{\mB}{\mathbb{B}}
\newcommand{\mC}{\mathbb{C}}
\newcommand{\mD}{\mathbb{D}}

\newcommand{\mN}{\mathbb{N}}

\newtheorem{theorem}{Theorem}[section]

\newtheorem{corollary}[theorem]{Corollary}

\newtheorem{observation}[theorem]{Observation}

\theoremstyle{definition}

\theoremstyle{definition}
\newtheorem{definition}[theorem]{Definition}

\theoremstyle{definition}

\begin{document}

\keywords{coherent ring, polydisc Wiener algebra, several complex variables}

\subjclass{Primary 32A38; Secondary 46J15, 46J20, 30H05, 13J99}

\title[Short  proof of the noncoherence of rings]{Noncoherence of some rings of holomorphic functions in several variables
as an easy consequence of the one-variable case}

\author{Raymond Mortini}
  \address{
 Universit\'{e} de Lorraine\\
 D\'{e}partement de Math\'{e}matiques et
Institut \'Elie Cartan de Lorraine,  UMR 7502\\
 Ile du Saulcy\\
 F-57045 Metz, France}
 \email{raymond.mortini@univ-lorraine.fr}

\author{Amol Sasane}
\address{Department of Mathematics, London School of Economics, Houghton Street, London WC2A 2AE, U.K.}
\email{sasane@lse.ac.uk}

\begin{abstract}
Using the facts that the disk algebra and the Wiener algebra are not coherent,
we prove that the polydisc algebra, the ball algebra and the Wiener algebra of the polydisc
are not coherent.
\end{abstract}

\maketitle

\section{Introduction}

Let us recall the notion of a coherent ring.

\begin{definition}
A unital commutative ring $R$ is said to be {\em coherent} if
 the intersection of any two finitely generated
ideals in $R$ is finitely generated, and for every $a\in R$, the annihilator
 $\textrm{Ann}(a):=\{x\in R: ax=0\}$ is finitely generated.
\end{definition}

We refer the reader
to the article \cite{Gla} for the relevance of the property
of coherence in commutative algebra. We are mainly interested in the question of coherence of an important nonuniform
algebra of holomorphic functions in several variables, namely the
Wiener algebra of the polydisc, defined below.

 \begin{definition}
  Let $\mD:=\{z\in \mC: |z|< 1\}$ and $\overline{\mD}:=\{z\in \mC: |z|\leq 1\}$.
  Let $n\in \mN$. The {\em Wiener algebra}  $W^+(\overline{\mD}^n)$ is the Banach algebra defined by
 $$
 W^+(\overline{\mD}^n)\!=\!\left\{ \! \sum_{k_1=0}^\infty\! \cdots \!\sum_{k_n=0}^\infty \! a_{(k_1,\cdots, k_n)} z_1^{k_1}\cdots z_n^{k_n}:
 \sum_{k_1=0}^\infty\! \cdots\! \sum_{k_n=0}^\infty\! |a_{(k_1,\cdots, k_n)}|<\infty\right\},
 $$
  with pointwise addition and multiplication, and the $\|\cdot\|_1$-norm given by
 $$
 \|f\|_{1}= \sum_{k_1=0}^\infty \!\cdots\! \sum_{k_n=0}^\infty\! |a_{(k_1,\cdots, k_n)}|, \;\;
 f=\displaystyle \!\sum_{k_1=0}^\infty \!\cdots \!\sum_{k_n=0}^\infty \!a_{(k_1,\cdots, k_n)} z_1^{k_1}\cdots z_n^{k_n}.
 $$
  The {\em polydisc algebra} $A(\overline{\mD}^n)$ is the Banach algebra
  of all continuous functions $f:\overline{\mD}^n\rightarrow \mC$ which are holomorphic in $\mD^n$,
  with pointwise addition and multiplication, and the supremum norm $\|\cdot\|_\infty$ given by
  $$
  \|f\|_\infty:=\sup_{(z_1,\cdots, z_n)\in \mD^n}|f(z_1,\cdots, z_n)|, \quad f\in A(\overline{\mD}^n).
  $$
 The ball algebra $A(\overline{\mB_n})$ is defined similarly, with the polydisc $\overline{\mD}^n$ replaced by the ball
  $
 \overline{\mB_n}:=\{(z_1\cdots,z_n)\in \mC^n :|z_1|^2+\cdots+|z_n|^2\leq 1\}.
 $
 \end{definition}

 McVoy  and Rubel  \cite{McVRub} showed that the disc algebra $A(\overline{\mD})$ is not coherent,
 but the question of coherence of the Wiener algebra was left open there.
 This was answered in
  \cite{Morvon}, where it was shown that the Wiener algebra $W^+(\overline{\mD})$ is not coherent. Using this, we extend that result
  to the Wiener algebra of the polydisc (which is new), and en route give
simplified proofs of the noncoherence of the polydisc algebra and the ball algebra (as opposed to the
rather technical proofs given in \cite{Ama} and \cite{Hic}).

 That the disc algebra is not coherent, can also be shown in a much easier way than
 in \cite{McVRub}. Let $I=(1-z)$ and $J=(1-z) S(z)$ where $S$ is given by $S(z)=\exp(-(1+z)/(1-z))$.
 Then $I\cap J=(1-z) S(z) M(1)$, where $M(1)$ is the maximal ideal of all those functions
 in $A(\mD)$ that vanish at $1$.  Since $M(1)=M(1)M(1)$ (just use that $f=\varphi F=
 (\varphi \sqrt F) \sqrt F$, where $\varphi F$ is the inner-outer factorization of $f$),
 $M(1)$ cannot be finitely generated by Nakayama's Lemma \cite[Theorem~76]{Kap}.

We also mention that the Wiener algebra $W^+(\overline{\mD}^n)$ arises naturally in Control Theory as a
class of stable transfer functions of discrete multidimensional systems, namely those linear time-invariant systems
which, for all $1\leq p\leq +\infty $, map inputs in $ \ell^p(\mN^n)$ to outputs in $ \ell^p(\mN^n)$.
 Noncoherence of the ring  $W^+(\overline{\mD}^n)$ then has consequences in the stabilization problem via the
 factorization approach; see \cite{Qua}.

\section{Coherence by reducing the number of variables}

Our main result is the following:

\begin{observation}
\label{main_theorem_A}
Let $\calR_n$ be the class of all commutative unital subrings of $\mC^{\mD^n}$ or of $\mC^{\mB_n}$ under the
usual pointwise operations. For $R_n\in\calR_n$ and $f\in R_n$,  set
 $
(D f )(z_1 ) := f (z_1 , 0, \cdots , 0)
$
and $D R_n := \{ D f : f \in R_n  \}$. For $f\in D R_n$, define $U f$ by $(U f )(z_1,\cdots, z_n):=f(z_1)$.
 If for all $f\in R_n$, $U D f\in R_n$, then the coherence of $R_n \in\calR_n$  implies the coherence of the ring $D R_n$.
\end{observation}
(The notation $D,U$ is chosen so as to indicate going down, respectively up, in the number of variables. It is clear that
for $f\in DR_n$, $DUf=f$.)

\begin{proof} Let $I$ and $J$ be two finitely generated ideals in $D R_n$.
 Let $I_n$ and $J_n$ be the ideals in $R_n$ generated by the functions $U f_k$ ($k=1,\cdots ,K $), respectively $U g_\ell$ ($\ell=1,\cdots ,L$),
 where the $f_1, \cdots, f_K $, respectively $g_1,\cdots ,g_L$, are the generators for $I$ and $J$. The coherence of
 $R_n$ implies that $I_n \cap  J_n$ is finitely generated. Let $\{p_1 , \cdots  , p_M\}$
 be a set of generators for the ideal $I_n \cap J_n$ in $R_n$. Then $\{D p_1 , \cdots , D p_M \}$
 is a set of generators for $I \cap  J$. Indeed, if $f\in I\cap J$, then in particular, $f\in I=(f_1,\cdots, f_K)$,
 so that there exist $\alpha_1,\cdots, \alpha_K\in DR_n$ such that $f=\alpha_1 f_1+\cdots+\alpha_K f_K$. Hence
 $$
 Uf=(U\alpha_1)(Uf_1)+\cdots +(U\alpha_K)(Uf_K)
 \in (Uf_1,\cdots, Uf_K)=I_n.
 $$
 Similarly, we see that $Uf\in J_n$ too. Hence $Uf\in I_n\cap J_n=(p_1,\cdots, p_M)$, and so there exist
 $\gamma_1,\cdots, \gamma_M\in R_n$ such that $Uf=\gamma_1 p_1+\cdots+ \gamma_M p_M$. Consequently,
  $$
  f=DUf = (D\gamma_1)(Dp_1)+\cdots+(D\gamma_M)(Dp_M)\in (D p_1 , \cdots , D p_M),
  $$
  and so $I\cap J\subset (D p_1 , \cdots , D p_M)$.
  Vice versa, for any $m=1,\cdots, M$, $p_m\in I_n$, in particular,
  $p_m\in I_n=(Uf_1,\cdots, Uf_K)$. So there exist $\theta_1,\cdots, \theta_K \in R_n$ such that
  $p_m=\theta_1(Uf_1)+\cdots+\theta_K(Uf_K)$. Hence
  \begin{eqnarray*}
   Dp_m&=&(D \theta_1)(DUf_1)+\cdots+ (D\theta_K)(DUf_K)\\
   &=& (D \theta_1)f_1+\cdots+ (D\theta_K)f_K
  \in (f_1,\cdots, f_K)=I.
  \end{eqnarray*}
  Similarly each $Dp_m$ also belongs to $J$. Hence $\{Dp_1,\cdots, Dp_M\}\subset I\cap J$, and so the ideal
   $(Dp_1,\cdots, Dp_M)\subset I\cap J$ too. It remains to verify the condition on the annihilators.

 Let $f\in D R_n$, and suppose that $g\in D R_n$ is such that $fg=0$. Then $U g \in \textrm{Ann}(U f )$. If $h_1,\cdots, h_r$ are generators for
 $ \textrm{Ann}(U f )$ in $R_n$, then $D h_1,\cdots, D h_r$ will be  generators for $\textrm{Ann}(f)$ in $D R_n$.
 Hence $D R_n$ is coherent.
\end{proof}

In light of this result, and the facts that $W^+(\overline{\mD})$ and $A(\overline{\mD})$ are not coherent, we obtain the
following consequences:

\begin{corollary}
\label{main_theorem_B}
For all $n\geq 1$, $W^+(\overline{\mD}^n)$, $A(\overline{\mD}^n)$ and $ A(\overline{\mB_n})$ are not coherent.
\end{corollary}

\end{document}